\documentclass[12pt]{article}
\usepackage{amssymb}
\usepackage{mathrsfs}
\usepackage{pst-poly}     
\usepackage{cite}
\usepackage{amsmath}
\usepackage {pstricks, pst-node}
\usepackage {graphics}

\newtheorem{thm}{Theorem}[section]
\newtheorem{lem}[thm]{Lemma}

\newtheorem{cor}[thm]{Corollary}

\newenvironment{pf}[1][Proof]{\noindent\textbf{#1.} }{\hfill\rule{1mm}{2mm}}

\makeatletter \@addtoreset{equation}{section} \makeatother

\begin{document}
\title{\bf Edge-Fault Tolerance of Hypercube-like Networks \thanks {The work was supported by
NNSF of China (No.11071233, 61272008).}}
\author
{Xiang-Jun Li$^{a,b}$
\quad Jun-Ming Xu$^a$
\footnote{Corresponding author: xujm@ustc.edu.cn (J.-M. Xu)}\\
{\small $^a$School of Mathematical Sciences}\\
{\small University of Science
and Technology of China,}\\
{\small  Wentsun Wu Key Laboratory of CAS, Hefei, 230026, China}  \\
{\small $^b$School of Information and Mathematics,}\\
{\small Yangtze University, Jingzhou, Hubei, 434023, China}\\
 }
\date{}
 \maketitle

\begin{abstract}

This paper considers a kind of generalized measure $\lambda_s^{(h)}$
of fault tolerance in a hypercube-like graph $G_n$ which contain
several well-known interconnection networks such as hypercubes,
varietal hypercubes, twisted cubes, crossed cubes and M\"obius
cubes, and proves $\lambda_s^{(h)}(G_n)= 2^h(n-h)$ for any $h$ with
$0\leqslant h\leqslant n-1$ by the induction on $n$ and a new
technique. This result shows that at least $2^h(n-h)$ edges of $G_n$
have to be removed to get a disconnected graph that contains no
vertices of degree less than $h$. Compared with previous results,
this result enhances fault-tolerant ability of the above-mentioned
networks theoretically.

\vskip6pt



\end{abstract}

\section{Introduction}
It is well known that interconnection networks play an important
role in parallel computing/communication systems. An interconnection
network can be modeled by a graph $G=(V, E)$, where $V$ is the set
of processors and $E$ is the set of communication links in the
network. For graph terminology and notation not defined here we
follow \cite{xu01}.

The connectivity of a graph $G$ is an important measurement for
fault-tolerance of the network, and the larger the connectivity is,
the more reliable the network is. Because the connectivity has some
shortcomings, Esfahanian~\cite{E89} proposed the concept of
restricted connectivity, Latifi {\it et al.}~\cite{lhp94}
generalized it to the restricted $h$-connectivity which can measure
fault tolerance of an interconnection network more accurately than
the classical connectivity. The concepts stated here are slightly
different from theirs.

For a given integer $h\,(\geqslant 0)$, an edge subset $F$  of a
connected graph $G$ is called an {\it $h$-super edge-cut}, or {\it
$h$-edge-cut} for short, if $G-F$ is disconnected and has the
minimum degree $\delta(G-F)\geqslant h$. The {\it $h$-super
edge-connectivity} of $G$, denoted by $\lambda^{(h)}_s(G)$, is
defined as the minimum cardinality over all $h$-edge-cuts of $G$. It
is clear that, for $h\geqslant 1$, if $\lambda^{(h)}_s(G)$ exists,
then $\lambda_s^{(h-1)}(G)\leqslant \lambda_s^{(h)}(G)$.

For any graph $G$ and a given integer $h$, determining
$\lambda_s^{(h)}(G)$ is quite difficult. In fact, the existence of
$\lambda_s^{(h)}(G)$ is an open problem so far when $h\geqslant 1$.
Only a little knowledge of results have been known on
$\lambda_s^{(h)}(G)$ for particular classes of graphs and small
$h$'s, such as, Xu~\cite{x00c} determined
$\lambda_s^{(h)}(Q_n)=2^h(n-h)$ for $h\leqslant n-1$.

It is widely known that the hypercube has been one of the most
popular interconnection networks for parallel computer/communication
system. However, the hypercube has the large diameter
correspondingly. To minimize diameter, various networks are proposed
by twisting some pairs of links in hypercubes, such as the varietal
hypercube $VQ_n$~\cite{cc94}, the twisted cube
$TQ_n$~\cite{AP89,AP91}, the locally twisted cube
$LTQ_n$~\cite{yem05a}, the crossed cube $CQ_n$~\cite{ef92,k97}, the
M\"obius cube $MQ_n$~\cite{cl95} and so on. Because of the lack of
the unified perspective on these variants, results of one topology
are hard to be extended to others. To make a unified study of these
variants, Vaidya {\it et al.}~\cite{vrs93} introduced the class of
hypercube-like graphs $HL_n$, which contains all the above-mentioned
networks. Thus, the hypercube-like graphs have received much
attention in recent years~\cite{ct07, cth03,lthh07, plk07, wl12,wwx10}.

In this paper, we determine $\lambda_s^{(h)}(G_n)= 2^h(n-h)$ for any
$G_n\in HL_n$ and $0\leqslant h\leqslant n-1$. Our result contains
many know conclusions and enhances the fault-tolerant ability of the
hypercube-like networks theoretically.

The proof of this result is in Section 3 by the induction on $n$ and
a new technique. Section 2 recalls the definition and Section 4
gives a conclusion on our work.

\section{Hypercube-like graphs}

Let $G_0=(V_0,E_0)$ and $G_1=(V_1,E_1)$ be two disjoint graphs with
the same order,  $\sigma$ a bijection from $V_0$ to $V_1$. A 1-1
connection between $G_0$ and $G_1$ is defined as an edge-set
$M_{\sigma}=\{x\sigma(x)|\ x\in V_0, \sigma(x)\in V_1\}$. And let
$G_0\oplus G_1$ denote $G=(V_0\cup V_1, E_0\cup E_1\cup
M_{\sigma})$. Clearly, $M_{\sigma}$ is a perfect matching of $G$.
Note that the operation $\oplus$ may generate different graphs
depending on the bijection $\sigma$.

Applying the above operation $\oplus$ repeatedly, a set of
$n$-dimensional {\it hypercube-like graphs}, denoted by $HL_n$, can
be recursively defined as follows.

(1) $HL_0=\{G_0\}$, where $G_0=K_1$, which is a single vertex;

(2) $G_n\in HL_{n}$ if and only if $G_n=G_{n-1}\oplus G'_{n-1}$ for
some $G_{n-1}, G'_{n-1}\in HL_{n-1}$.

\vskip6pt

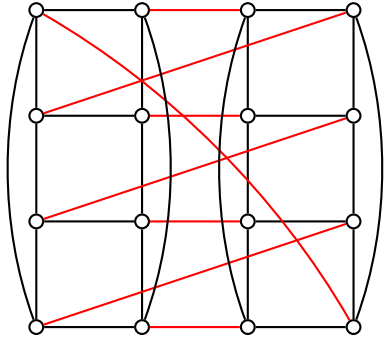
\begin{figure}[ht]
\begin{center}
\psset{unit=20pt}
\begin{pspicture}(-4,-3.5)(4,4)

\cnode(-3,3){3pt}{a1} \cnode(-1,3){3pt}{a2} \ncline{a1}{a2}
\cnode(-1,1){3pt}{a3} \cnode(-3,1){3pt}{a4} \ncline{a2}{a3}\ncline{a3}{a4}\ncline{a1}{a4}

\cnode(1,3){3pt}{b1} \cnode(3,3){3pt}{b2} \ncline{b1}{b2}
\cnode(3,1){3pt}{b3} \cnode(1,1){3pt}{b4} \ncline{b2}{b3}\ncline{b3}{b4}\ncline{b1}{b4}

\cnode(-3,-1){3pt}{c1} \cnode(-1,-1){3pt}{c2} \ncline{c1}{c2}
\cnode(-1,-3){3pt}{c3} \cnode(-3,-3){3pt}{c4} \ncline{c2}{c3}\ncline{c3}{c4}\ncline{c1}{c4}

\cnode(1,-1){3pt}{d1} \cnode(3,-1){3pt}{d2} \ncline{d1}{d2}
\cnode(3,-3){3pt}{d3} \cnode(1,-3){3pt}{d4} \ncline{d2}{d3}\ncline{d3}{d4}\ncline{d1}{d4}

\ncline{a4}{c1} \ncline{a3}{c2} \ncline{b4}{d1}\ncline{b3}{d2}
\ncline[linecolor=red]{a2}{b1} \ncline[linecolor=red]{a3}{b4}
\ncline[linecolor=red] {c2}{d1} \ncline[linecolor=red]{c3}{d4}
\ncline[linecolor=red]{b2}{a4} \ncline[linecolor=red]{b3}{c1}
\ncline[linecolor=red]{d2}{c4}
\nccurve[angleA=-30,angleB=120,linecolor=red]{a1}{d3}
\nccurve[angleA=-110,angleB=110]{a1}{c4}
\nccurve[angleA=-70,angleB=70]{a2}{c3}
\nccurve[angleA=-110,angleB=110]{b1}{d4}
\nccurve[angleA=-70,angleB=70]{b2}{d3}

\end{pspicture}
\caption{\label{f1}\footnotesize {A hypercube-like graph in $HL_4$}}
\end{center}
\vskip-20pt
\end{figure}

It is clear that for a graph $G_n\in HL_n$, $G_n$ is an $n$-regular connected graph
of order $2^n$. A hypercube-like graph in $HL_4$ is shown in Fig~\ref{f1}.

By definitions, the hypercube $Q_n=Q_{n-1}\oplus Q_{n-1}$, the
varietal hypercube $VQ_n=VQ_{n-1}\oplus VQ_{n-1}$, the twisted cube
$TQ_n=TQ_{n-1}\oplus TQ_{n-1}$, the locally twisted cube
$LTQ_n=LTQ_{n-1}\oplus LTQ_{n-1}$, the crossed cube
$CQ_n=CQ_{n-1}\oplus CQ_{n-1}$, the M\"obius cube
$MQ_n=MQ_{n-1}\oplus MQ_{n-1}$. Thus, $\{Q_n, VQ_n, TQ_n, LTQ_n,
CQ_n, MQ_n\}\subseteq HL_n$.

For convenience, for a graph $G$, we write $|G|$ for $|V(G)|$, for a
subgraph $X\subseteq G$, write $X$ for $V(X)$. For $G_n\in HL_n$, we
write $M_n$ for $M_{\sigma}$. Let $I_n=\{0,1,\ldots,n\}$. From the
definition of $HL_n$, it is easy to see that, for two given integers
$n$ and $h\in I_{n}$ and a given graph $G_n\in HL_n$, there is a
graph $G_h\in HL_h$ from which $G_n$ can be obtained by repeating
$n-h$ times of the operation $\oplus$, and denote $G_n=G_h^{n-h}$.
Moreover, for any edge $e$ in $G_n$, there exists some $h\in I_n$
such that $e$ is an edge in $G_h\in HL_h$ and $e\in M_h$.

\section{Main results}

In this section, our aim is to prove that
$\lambda_s^{(h)}(G)=2^h(n-h)$ for any $h\in I_{n-1}$.

\begin{lem}\label{lem3.1}
$\lambda_s^{(h)}(G_n)\leqslant 2^h(n-h)$ for any $G_n\in HL_n$ and
$h\in I_{n-1}$.
\end{lem}

\begin{pf}
Let $G_n\in HL_n$. Then there exist a graph $G_h\in HL_h$ such that
$G_n=G_h^{n-h}$. Let $F$ be the set of edges between $G_h$ and
$G_n-G_h$. Then $F$ is an edge-cut of $G_n$. Since $G_n$ is
$n$-regular and $G_h\subseteq G_n$ is $h$-regular, $|F|=|G_h|
(n-h)=2^h(n-h)$.

We now show that $F$ is an $h$-edge-cut of $G_n$ by proving
$\delta(G_n-F)\geqslant h$. Let $x$ be a vertex in $G_n-F$. If
$x$ is in $G_h$, then $x$ in $G_n-F$ has degree $h$ clearly since
$G_h\in HL_h$. If $x$ is in $G_n-G_h$, it can be matched at most one
vertex in $G_h$ by the matching $M_i$ for some $i\in
\{h+1,h+2,\ldots,n\}$, which implies that $x$ has degree at least
$(n-1\geqslant)\,h$ in $G_n-F$. By the arbitrariness of $x$,
$\delta(G_n-F)\geqslant h$, which shows that $F$ is an $h$-edge-cut
of $G$, and so
 $$
 \lambda_s^{(h)}(G)\leqslant |F|=2^h(n-h)
 $$
The lemma follows.
\end{pf}

\vskip 8pt For $G_n\in HL_n$, let $X \subseteq G_n$ be a subgraph of
$G_n$, $Y=G_n-X$. $E_n(X)$ denote the set of edges between $X$ and
$Y$ in $G_n$. For convenience, let $G_n=H_0\oplus H_1$, where
$H_0=G_{n-1}$ and $H_1=G'_{n-1}$. For $i\in I_1$, let
 $$
 X_i=X\cap H_i,\ \  Y_i=Y\cap H_i,\ \ F_i=E_n(X)\cap H_i\ \ {\rm and}\ F_2=E_n(X)\cap M_n.
 $$

\begin{lem}\label{lem3.2}\
$|X|\geqslant 2^h$ if $\delta(X)\geqslant h$ for $h\in I_{n}$.
\end{lem}

\begin{pf}
We proceed by induction on $n\,(\geqslant 1)$ for fixed $h\in
I_{n-1}$ by the recursive structure of $G_n$. Clearly, if $n=1$ then
the conclusion is true for any $h\in I_1$. Assume that the
conclusion holds for $n-1$ with $n\geqslant 2$.

If $X\subseteq H_0$ or $X\subseteq H_1$, then we have done by our
hypothesis. Assume $X_i=X\cap H_i\ne\emptyset$ for each $i\in I_1$
below. Then $\delta(X_i)\geqslant h-1$ in $H_i$ for each $i\in I_1$
since $\delta(X)\geqslant h$ and there is at most one edge linking a
vertex in $X_0$ and a vertex in $X_1$ in $G_n$. By our hypothesis,
$|X_i|\geqslant 2^{h-1}$ for each $i\in I_1$. It follows that
 $$
 |X|=|X_0|+|X_1|\geqslant 2\cdot 2^{h-1}=2^{h}.
 $$
By the induction principle, the lemma follows.
\end{pf}

\begin{lem}\label{lem3.5}\
$|X|+|E_n(X)| \geqslant 2^{h}(n+1-h)$ if $\delta(X)\geqslant h$ for any
$h\in I_{n-1}$.
\end{lem}

\begin{pf}
We proceed by induction on $n\,(\geqslant 1)$ for fixed $h\in
I_{n-1}$. Clearly, the conclusion hold for $n=1$. Assume the
induction hypothesis for $n-1$. There are two cases.

\vskip6pt

{\bf Case 1.} $X\subseteq H_0$ or $X\subseteq H_1$.

We assume  $X\subseteq H_0$ without loss of generality. If $X=H_0$,
then $h=n-1$, and $|E_n(X)|=|H_0|=2^{n-1}$, so the conclusion is hold. Assume
$X\subset H_0$ below. Then $h\leqslant n-2$. Since every vertex in
$X$ has exactly one neighbor in $H_1$ matched by a perfect matching
$M_n$, we have $|F_2|=|X|$, and so $|F_2|=|X|\geqslant 2^h$ by
Lemma~\ref{lem3.2}.

Since $\delta(X)\geqslant h$ and $h\in I_{n-2}$, using the induction
hypothesis in $H_0$, we have $|X|+|E_{n-1}(X)|\geqslant 2^h(n-h)$.
Combining this with $|F_2|\geqslant 2^h$, we have
 $$
\begin{array}{rl}
|X|+|E_n(X)|&=|X|+|E_{n-1}(X)|+|F_2|\\
&\geqslant 2^h(n-h)+2^h\\
&=2^h(n+1-h).
 \end{array}
 $$
and so the conclusion holds.

 \vskip6pt

 {\bf Case 2.} $X_i=X\cap H_i\ne\emptyset$ for each $i\in I_1$.

Since $\delta(X)\geqslant h$ in $G_n$, $\delta(X_i)\geqslant h-1$ in
$H_i$ for $i\in\{0,1\}$. Using the induction hypothesis in $H_i$ for
$i\in\{0,1\}$, we have
 $$
\begin{array}{l}
|X_i|+|E_{n-1}(X_i)| \geqslant 2^{h-1}(n +1-h)\ \ {\rm for\ each}\
i\in I_1.
\end{array}
 $$
It follows that
  $$
\begin{array}{l}
|X|+|E_{n}(X)|\geqslant
|X_0|+|E_{n-1}(X_0)|+|X_1|+|E_{n-1}(X_1)|\geqslant 2^{h}(n+1-h).
\end{array}
 $$

By the induction principle, the lemma follows.
\end{pf}

\begin{lem}\label{lem3.7}\
$|E_n(X)| \geqslant 2^{h}(n-h)$ if $\delta(X)\geqslant h$ and
$\delta(Y)\geqslant h$ for any $h\in I_{n-1}$.
\end{lem}

\begin{pf}
By symmetry of $X$ and $Y$, we can assume $|X|\leqslant |Y|$. We
prove the conclusion by induction on $n(\geqslant 1)$. The
conclusion is true for $n=1$ clearly. Assume the induction
hypothesis for $n-1$. There are two cases.

\vskip6pt

{\bf Case 1.} $X\subseteq H_0$ or $X\subseteq H_1$.

Assume  $X\subseteq H_0$ without loss of generality. If $X=H_0$,
then $h=n-1$, and so the conclusion is straightforward. Assume
$X\subset H_0$ below. Then $h\leqslant n-2$.

Clearly, $|F_2|= |X|$ since every vertex in $X$ has exactly one
neighbor in $H_1$ matched by a perfect matching $M_n$. Since
$\delta(X)\geqslant h$, $X \subset H_0$ and $h\in I_{n-2}$, using
Lemma~\ref{lem3.5} in $H_0$, we have
$$
\begin{array}{rl}
|X|+|E_{n-1}(X)| \geqslant 2^{h}(n-h)
\end{array}
$$
It follows that
$$
\begin{array}{rl}
|E_n(X)|&=|E_{n-1}(X)|+|F_2|\\
&= |E_{n-1}(X)|+|X|\\
&\geqslant 2^h(n-h),
 \end{array}
$$
and so the conclusion holds.
 \vskip6pt

 {\bf Case 2.}  $X_i=X\cap H_i\ne\emptyset$ for each $i\in I_1$.

Clearly, $Y_0 \ne \emptyset$ and $Y_1\ne \emptyset$ since
$|X|\leqslant |Y|$, and $\delta(X_i)\geqslant h-1$ and
$\delta(Y_i)\geqslant h-1$ in $H_i$ for $i \in \{0,1\}$ since
$\delta(X)\geqslant h$ and $\delta(Y)\geqslant h$ in $G_n$. Then
$X_i\subseteq H_i$ and $Y_i\subseteq H_i$, satisfy our
hypothesis for $i \in \{0,1\}$. By the induction hypothesis, we have
 $$
\begin{array}{l}
|F_i|=|E_{n-1}(X_i)|\geqslant 2^{h-1}(n-h)\ \ {\rm for\ each}\
i\in \{0,1\},
\end{array}
 $$
It follows that
 $$
 |E_n(X)| \geqslant |F_0|+|F_1|\geqslant 2\cdot 2^{h-1}(n-h)=2^{h}(n-h),
 $$
and so the conclusion holds.

By the induction principle, the lemma follows.
\end{pf}

\begin{thm}\label{thm3.8}
 $\lambda_s^{(h)}(G_n)= 2^h(n-h)$ for any $G_n\in HL_n$ and any $h\in I_{n-1}$.
\end{thm}

\begin{pf}
Let $G_n\in HL_n$. By Lemma~\ref{lem3.1}, we need only to show that
$\lambda_s^{(h)}(G_n) \geqslant 2^h(n-h)$ for any $h\in I_{n-1}$.
Let $F$ be an $h$-edge-cut of $G_n$ with $|F|=\lambda_s^{(h)}(G_n)$,
$X$ a connected component of $G_n-F$, and $Y=G_n-X$. Clearly,
$\delta(X)\geqslant h$ and $\delta(Y)\geqslant h$ since $F$ is an
$h$-edge-cut. By Lemma~\ref{lem3.7}, we immediately have
 $$
 \lambda_s^{(h)}(G_n)=|F|=|E_n(X)|\geqslant 2^h(n-h),
 $$
and so the theorem follows.
\end{pf}

\begin{cor}
If $G_n\in \{Q_n, VQ_n, CQ_n, MQ_n, TQ_n, LTQ_n\}$, then
$\lambda_s^{(h)}(G_n)=2^h(n-h)$ for any $h\in I_{n-1}$.
\end{cor}

\section{Conclusions}

In this paper, we consider the generalized measures of edge fault
tolerance for the hypercube-like networks, called the $h$-super
edge-connectivity $\lambda_s^h$. For the hypercube-like graph $G\in
HL_n$, we prove that $\lambda_s^{(h)}(G)= 2^h(n-h)$ for any $h\in
I_{n-1}$. The results show that at least $2^h(n-h)$ edges of $G$
have to be removed to get a disconnected graph that contains no
vertices of degree less than $h$. Thus, when the hypercube-like
networks is used to model the topological structure of a large-scale
parallel processing system, these results can provide more accurate
measurements for fault tolerance of the system.

Similarly, we can define $h$-super connectivity $\kappa_s^h(G)$ of a
connected graph $G$ by considering vertices rather than edges. One
may ask if $\kappa_s^{(h)}(G)= 2^h(n-h)$ for any $h\in I_{n-1}$ with
$0\leqslant h\leqslant n-1$, or how many vertices of $G\in HL_n$
have to be removed to get a disconnected graph that contains no
vertices of degree less than $h$. In fact, there is some graph $G\in
HL_n$ such that $\kappa_s^{(h)}(G)$ does not exist, that is, no
matter how we remove the vertices, we can not get a disconnected
graph that contains no vertices of degree less than $h$. The graph
shown in Fig~\ref{f1} is an example for $h=2$. It is worth while to
research the existence of $\kappa_s^{(h)}(G)$ for some $G\in HL_n$
or $h\in I_{n-1}$, and to determine $\kappa_s^{(h)}(G)$ if
$\kappa_s^{(h)}(G)$ exists.

\end{document}